\documentclass{amsart}

\usepackage{amssymb}
\usepackage{graphicx}

\newtheorem{theorem}{Theorem}[section]
\newtheorem{lemma}[theorem]{Lemma}
\newtheorem{proposition}[theorem]{Proposition}

\newtheorem{corollary}[theorem]{Corollary}

\newtheorem{_algorithm}[theorem]{Algorithm}

\newtheorem{_definition}[theorem]{Definition}
\newenvironment{definition}{\begin{_definition}\rm}{\end{_definition}}

\newtheorem{_remark}[theorem]{\it Remark}
\newenvironment{remark}{\begin{_remark}\rm}{\end{_remark}}

\newtheorem{_example}[theorem]{Example}
\newenvironment{example}{\begin{_example}\rm}{\end{_example}}

\newtheorem{_assumption}[theorem]{Assumption}

\newtheorem{_construction}[theorem]{Construction}

\newtheorem{_claim}[theorem]{Claim}

\newtheorem{_conjecture}[theorem]{Conjecture}

\newtheorem{_problem}[theorem]{Problem}

\numberwithin{equation}{section}
\numberwithin{table}{section}
\numberwithin{figure}{section}

\newcommand{\C}{\mathord{\mathbb C}}

\renewcommand{\P}{\mathord{\mathbb  P}}
\newcommand{\Q}{\mathord{\mathbb  Q}}
\newcommand{\R}{\mathord{\mathbb R}}
\newcommand{\Z}{\mathord{\mathbb Z}}

\newcommand{\CCC}{\mathord{\mathcal C}}

\newcommand{\HHH}{\mathord{\mathcal H}}

\newcommand{\LLL}{\mathord{\mathcal L}}
\newcommand{\MMM}{\mathord{\mathcal M}}

\newcommand{\OOO}{\mathord{\mathcal O}}

\newcommand{\QQQ}{\mathord{\mathcal Q}}

\newcommand{\maprightsp}[1]{\; \smash{\mathop{\; \longrightarrow \; }\limits\sp{#1}}\; }

\newcommand{\mapdownsurj}{
\hbox{$\bigm\downarrow$}
\llap{\hbox{\raise 2pt\hbox{$\bigm\downarrow$}}}%
\vstrechmapdown
}

\newcommand{\inj}{\hookrightarrow}

\newcommand{\surj}{\mathbin{\to \hskip -7pt \to}}

\newcommand{\isom}{\mathbin{\,\raise -.6pt\rlap{$\to$}\raise 3.5pt%
\hbox{\hskip .3pt$\mathord{\sim}$}\,}}

\newcommand{\set}[2]{\{\; {#1} \; \mid \; {#2} \;  \}}

\newcommand{\map}[3]{ #1 \, : \, #2 \, \to \, #3}
\newcommand{\mapisom}[3]{ #1 \, : \, #2 \; \isom \; #3}

\newcommand{\shortmap}[3]{ #1  :  #2 \to #3}

\newcommand{\sm}{\setminus}
\newcommand{\st}{\subset}

\newcommand{\sprime}{\sp\prime}

\newcommand{\sptimes}{\sp{\times}}

\newcommand{\dual}{\sp{\vee}}

\newcommand{\inv}{\sp{-1}}

\newcommand{\Ker}{\operatorname{\rm Ker}\nolimits}

\renewcommand{\Im}{\operatorname{\rm Im}\nolimits}

\newcommand{\Aut}{\operatorname{\rm Aut}\nolimits}

\newcommand{\pr}{\operatorname{\rm pr}\nolimits}
\newcommand{\Sing}{\operatorname{\rm Sing}\nolimits}

\newcommand{\rank}{\operatorname{\rm rank}\nolimits}

\newcommand{\Pt}{\P^2}
\newcommand{\pione}{\pi_1}

\newcommand{\Emb}{\operatorname{\rm Emb}\nolimits} 
\newcommand{\Pic}{\operatorname{\rm Pic}\nolimits}

\newcommand{\GL}{\operatorname{\it{GL}}\nolimits}
\newcommand{\SL}{\operatorname{\it{SL}}\nolimits}

\newcommand{\wt}[1]{\widetilde{#1}}

\newcommand{\ang}[1]{\langle #1\rangle}

\newcommand{\rmand}{\textrm{and}}

\newcommand{\quand}{\quad\rmand\quad}

\newenvironment{rmenumerate}
{\begin{enumerate}

}
{\end{enumerate}}

\newcommand{\mystruth}[1]{\phantom{\hbox{\vrule height #1}}}
\newcommand{\mystrutd}[1]{\phantom{\hbox{\vrule depth #1}}}

\newcommand{\erase}[1]{}

\newcommand{\Cl}{\operatorname{\it Cl}\nolimits} 

\newcommand{\Cinf}{\CCC^\infty}

\newcommand{\tor}{_{\mathord{\rm{tor}}}}
\newcommand{\tf}{^{\mathord{\rm{tf}}}}

\newcommand{\NC}{N_C}
\newcommand{\MC}{M_C}
\newcommand{\BU}{B_U}

\newcommand{\NCspsim}{\widetilde{N}_C}
\newcommand{\MCspsim}{\widetilde{M}_C}
\newcommand{\BUspsim}{\widetilde{B}_U}

\newcommand{\spnu}{^{(\nu)}}

\newcommand{\lat}{\Lambda}
\newcommand{\tQ}{\otimes \Q}
\newcommand{\tC}{\otimes \C}
\newcommand{\rootlat}{\Sigma}

\newcommand{\relQ}{\mathord{\sim}}

\newcommand{\Ms}{\mathord{M\hskip -1pt s}}
\newcommand{\Ns}{\mathord{N\hskip -1pt s}}
\newcommand{\Ls}{\mathord{Ls}}

\newcommand{\cpersp}{\mathord{c\hskip.7pt\Omega}}
\newcommand{\cpersps}{\mathord{c\hskip.7pt\Omega s}}

\newcommand{\polar}{\LLL}
\newcommand{\sbXL}{_{(X, \polar)}}

\newcommand{\Orb}{U}

\newcommand{\Mat}{Q}
\newcommand{\Mats}{\QQQ}
\newcommand{\Lat}{\Lambda}
\newcommand{\oriLat}{\widetilde{\Lambda}}

\newcommand{\AutFh}{O_{F, h} (M^0)}
\newcommand{\AutFhM}{O_{F, h, M} (M^0)}

\newcommand{\spsharp}{\sp{\raise -1pt \hbox{${}_\sharp$}}}

\newcommand{\Km}{\mathord{\rm Km}}

\begin{document}

\title[Arithmetic Zariski pairs]{On arithmetic Zariski pairs in degree $6$}

\author{Ichiro Shimada}
\address{
Department of Mathematics,
Faculty of Science,
Hokkaido University,
Sapporo 060-0810,
JAPAN
}
\email{shimada@math.sci.hokudai.ac.jp
}

\dedicatory{Dedicated to Professor Mutsuo Oka for his sixtieth birthday.}


\begin{abstract}
We define a topological invariant of  complex projective plane curves.
As an application,
we present  new examples of arithmetic Zariski pairs.
\end{abstract}

\maketitle

\section{Introduction}\label{sec:Introduction}
In this paper, we mean by a \emph{plane curve}
a complex reduced (possibly reducible) projective plane curve.
The following definition is due to Artal-Bartolo~\cite{MR1257321}:
\begin{definition}
A pair $(C, C\sprime)$ of plane curves 
of the same degree is called a \emph{Zariski pair}
if there exist  tubular neighborhoods $T\subset \Pt$ of $C$ and $T\sprime\subset \Pt$ of $C\sprime$
such that $(T, C)$ and $(T\sprime, C\sprime)$ are diffeomorphic,
while $(\Pt, C)$ and $(\Pt, C\sprime)$ are not homeomorphic.
\end{definition}
The first example of Zariski pairs was 
studied by Zariski~\cite{MR1506719}
in order to show that 
an equisingular family of plane curves need not be connected.
Zariski considered a  six-cuspidal sextic curve  $C$ with
the six  cusps lying on a conic,
and proved that $\pione (\Pt\sm C)$ is isomorphic to
the free  product  of   cyclic groups of order $2$ and $3$.
He then 
showed that,  if   there exists a six-cuspidal sextic curve  $C\sprime$  with
the six  cusps \emph{not} lying on a conic,
then $\pione (\Pt\sm C\sprime)$ is not isomorphic to $\pione (\Pt\sm C)$.
Oka~\cite{MR1167373} completed Zariski's work
by constructing explicitly a non-conical six-cuspidal sextic curve $C\sprime$,
and showed that $\pione (\Pt\sm C\sprime)$
is a cyclic group of order $6$.
Therefore 
the moduli space $\MMM(6A_2)$ of
plane sextics possessing  six  cusps as their only singularities 
 has at least two connected components
that are distinguished by the fundamental groups of the complements.
(See~\cite{MR1421396} for a simple construction of the pair $(C, C\sprime)$.)
Recently, Degtyarev~\cite{degtyarev-2005} showed that 
$\MMM(6A_2)$ has exactly two connected components.
\par
\medskip
Many examples of Zariski pairs have been known now.
The standard method to distinguish $(\Pt, C)$ and $(\Pt, C\sprime)$ topologically
is to compare the fundamental groups
of the complements.
The fundamental groups are calculated directly
by Zariski-van Kampen theorem (see~\cite{MR1952329}),
or they are compared by some indirect methods;
for example, 
by means of  Alexander polynomials,
or by proving  (non-)existence of finite \'etale Galois coverings of the complements
with a given Galois group.
\begin{definition}
Plane curves $C$ and $C\sprime$ are said to be \emph{conjugate}
if there exist a homogeneous polynomial $\Phi(x_0, x_1, x_2)$ 
of complex coefficients and an automorphism 
$\sigma$ of the field $\C$ such that we have 
$$
C=\{\Phi=0\}\quand C\sprime=\{\Phi^\sigma=0\},
$$
where $\Phi^\sigma$ is 
the polynomial obtained from $\Phi$ by applying $\sigma$ to the coefficients of $\Phi$.
\end{definition}
\begin{definition}
A Zariski pair  $(C, C\sprime)$ is called  an \emph{arithmetic Zariski pair}
if $C$ and $C\sprime$ are conjugate.
\end{definition}
A difficulty in constructing   examples of arithmetic Zariski pairs comes from  the fact that, 
if $C$ and $C\sprime$ are conjugate,
then $\pione (\Pt\sm C)$ and $\pione (\Pt\sm C\sprime)$ have the same pro-finite completions.
\par
\medskip
Artal~Bartolo, Carmona~Ruber, and Cogolludo~Agust{\'{\i}}n~(\cite{MR1980995}, ~\cite{MR2247887})
constructed an arithmetic Zariski pair in degree $12$.
They distinguished $(\Pt, C)$ and $(\Pt, C\sprime)$
by means of the \emph{braid monodromy}.
\par
\medskip
In this paper,
we introduce an invariant $\NC$ of the homeomorphism type of $(\Pt, C)$
for plane curves $C$ of even degree.
By means of this invariant,
we present some examples  of arithmetic  Zariski pairs in degree $6$.
\par
\medskip
In order to explain our examples more precisely, we introduce some terminologies.
A \emph{Dynkin type} is a finite formal sum
$$
R=\sum_{l\ge 1} a_l A_l + \sum_{m\ge 4} d_m D_m +\sum_{n=6}^{8} e_n E_n, 
$$
where $a_l, d_m$ and $e_n$ are non-negative integers,
almost all of which are zero.
The \emph{rank} of the Dynkin type $R$ is defined by
$$
\rank (R)=\sum a_l l + \sum d_m m +\sum e_n n.
$$
An \emph{$ADE$-sextic} is a plane curve of degree
$6$ with only simple singularities.
The \emph{type $R$ of an $ADE$-sextic $C$}
is the Dynkin type of the singularities of $C$.
Then $\rank (R)$ is equal to the total Milnor number of $C$,
and hence it is at most $19$.
We say that an $ADE$-sextic is 
a \emph{maximizing sextic} if the total Milnor number is $19$
(see Persson~\cite{MR805337}).
If $C$ is an $ADE$-sextic,
then the minimal resolution $X_C$ of the double covering   $Y_C\to \Pt$ that branches  exactly along $C$
is a  $K3$ surface.
When $C$ is a maximizing sextic,
our invariant $\NC$ of $(\Pt, C)$ coincides with 
the transcendental lattice of $X_C$.
 \par
 \medskip
Combining our main result
with the results
 of Artal-Bartolo,
  Carmona-Ruber,
  Cogolludo-Agust{\'{\i}}n~\cite{MR1900779},
 Degtyarev~\cite{degtyarev-2005} and 
 the  result in~\cite{tssK3},
 we obtain the following:
 \begin{theorem}
 There exists an arithmetic Zariski pair  of maximizing sextics 
  for each of the following  Dynkin types:
 $$
 {\rm (i)} \; A_{16}+A_{2}+A_1, 
 \quad
 {\rm (ii)} \; A_{16}+A_{3},
 \quad
 {\rm (iii)} \; A_{18}+A_{1},
 \quad
 {\rm (iv)} \; A_{10}+A_{9}.
 $$
  \end{theorem}
The plan of this paper is as follows.
In~\S\ref{sec:N},
we define an invariant $\NC$ for  curves $C$ on a smooth projective surface $S$
satisfying  certain conditions, 
and show that $\NC$ is in fact an invariant of 
the \emph{$\varGamma$-equivalence class} of $(S, C)$
(see Definition~\ref{def:Gamma}).
The $\varGamma$-equivalence
is an equivalence relation coarser than
the homeomorphism type
of $(S, C)$,
and finer 
than
the homeomorphism type
of $S\sm C$.
Applying the main result (Theorem~\ref{thm:main}) to the case  $S=\Pt$,
we obtain an invariant of
the homeomorphism type of  $(\Pt, C)$
for plane curves $C$ of  even degree.
In~\S\ref{sec:conn},
we review the theory of Degtyarev~\cite{degtyarev-2005}
on the connected components of 
the moduli space $\MMM(R)$ of $ADE$-sextics
with a given Dynkin type $R$.
In~\S\ref{sec:examples},
we calculate the connected components of $\MMM(R)$
for some $R$ with  $\rank (R)=19$.
Combining this calculation  with the result of~\cite[Theorem~5.8]{MR1900779}
and using our  invariant, 
we show that some pairs of conjugate maximizing sextics 
obtained in~\cite{MR1900779}
yield  examples  of arithmetic Zariski pairs
(the examples (i)-(iii) above).
In \S\ref{sec:singK3}, 
 we present another example 
of arithmetic Zariski pairs
constructed 
by means of 
the theory of Hilbert class fields
of  imaginary quadratic fields (the example (iv) above).
\par
\medskip
The first example of non-homeomorphic conjugate complex varieties was given by Serre~\cite{MR0166197}.
Since then, only few examples seem to have been treated~(e.g., Abelson~\cite{MR0349679}).
The argument of this paper provides us with
a new method to construct examples of non-trivial effects
of $\Aut (\C)$ on the topology of complex varieties.
\section{The  invariant $\NC$}\label{sec:N}
First we fix some notation and terminologies.
\par
\medskip
Let $A$ be a finitely generated $\Z$-module.
We denote by $A\tor$ the torsion subgroup of $A$,
and  by  $A\tf:=A/A\tor$
the \emph{torsion-free quotient} of $A$.
If
$\shortmap{b}{A\times A}{\Z}$
is a symmetric bilinear form on $A$, 
then $b$ induces
a symmetric bilinear form on $A\tf$
in the natural way.
\par
\smallskip

Let $A$ be a free $\Z$-module of finite rank,
and $A\sprime$ a submodule of $A$.
The \emph{primitive closure} of $A\sprime$ in $A$ is 
defined to be the intersection of 
$A\sprime\otimes \Q$ and  $A$ in $A\otimes \Q$.
We say that $A\sprime$ is \emph{primitive in $A$}
if the primitive closure of $A\sprime$ is equal to $A\sprime$.
\par
\smallskip
A \emph{lattice} is a free $\Z$-module $A$ of finite rank
equipped with a \emph{non-degenerate} 
symmetric bilinear form $A\times A\to \Z$.
Two lattices $A$ and $A\sprime$ are \emph{isomorphic} if 
there exists an isomorphism $A\isom A\sprime$ of $\Z$-modules
that preserves the symmetric bilinear forms.
The automorphism group of a lattice $A$ is denoted by $O(A)$.
If $A$ and $A\sprime$ are lattices,
then $A\perp A\sprime$ denotes the orthogonal direct-sum of $A$ and $A\sprime$.
\par
\smallskip
For a topological space $Z$, we denote by $H_2(Z)$ the homology group $H_2(Z,\Z)$.
When $Z$ is  an oriented  $\Cinf$-manifold with $\dim_{\R} (Z)=4$, we have 
the intersection pairing
$\shortmap{b_Z}{H_2 (Z)\times H_2(Z)}{\Z}$.
If we further assume that $Z$ is   compact,
then $H_2(Z)\tf$ becomes a lattice by $b_Z$.
\par
\medskip
Let $S$ be a smooth  complex projective surface such that 
 $$
 \Pic (S)\cong \Z
 \quand 
 \pione (S)=\{1\},
 $$
and  let $\HHH$ be  the line bundle on $S$ such that its class is the positive generator of $\Pic (S)$.
Let $d$ be a positive  even integer, and put
$$
\LLL:=\HHH\sp{\otimes d}
\quand
\MMM:=\HHH\sp{\otimes d/2}.
$$
An \emph{$\LLL$-curve} is a reduced (possibly reducible)
member of the complete linear system $|\LLL|$.
Let $C$ be an $\LLL$-curve
given by $s=0$, where $s$ is  a global section of $\LLL$, 
and let
$$
\map{\pi}{Y}{S}
$$
be the finite double covering
that branches exactly along $C$,
where $Y$ is the pull-back of the image of the global section $s$ 
by the squaring morphism  $\MMM\to\MMM\sp{\otimes 2}=\LLL$ over $S$.
Note that $Y$ is normal,
because $Y$ is a hypersurface in the total space of the line bundle  $\MMM$
with only isolated singular points~(Altman and Kleiman~\cite[Chapter VII, Corollary (2.14)]{MR0274461}).
Let
$$
\map{\rho}{X}{Y}
$$
be a proper birational morphism
from a smooth surface $X$ that induces an isomorphism $\rho\inv (Y\sm \pi\inv (C))\cong Y\sm \pi\inv (C)$.
We put
$$
\map{\phi:=\pi\circ \rho\;}{\;X\;}{\;S}.
$$
Then $\phi$
is an  \'etale double covering over $S\setminus C$.
We denote by
$$
\MCspsim \;\;\subset\;\; H_2(X)
$$
the submodule generated by the homology classes of
the integral components of $\phi\inv (C)\subset X$.
We then put
\begin{eqnarray*}
\NCspsim  &:=& \set{x\in H_2 (X)}{b_X(x, y)=0\;\;\textrm{for any}\;\; y \in \MCspsim } \quand\\
\NC &:=& (\NCspsim  )\tf\;\;\subset\;\; H_2(X)\tf.
\end{eqnarray*}
Note that $\NC$ is primitive in $H_2(X)\tf$.
\begin{lemma}\label{lem:MN}
The restriction of $b_X$ to   $\NC$ is non-degenerate.
\end{lemma}
\begin{proof}
Since $H_2(X)\tf$ is a lattice by $b_X$,
and  $N_C$ is the orthogonal complement of $(\MCspsim )\tf$ in $H_2(X)\tf$,
it is enough to show that the restriction of $b_X$ to  $(\MCspsim )\tf$
is non-degenerate.
Let $h\in H_2(X)$ be the first Chern class of the line bundle $\phi^*(\HHH)$.
(We have a canonical  isomorphism $H_2(X)\cong H^2(X)$.)
Since $b_X(h, h)>0$,
the $\Z$-module $\ang{h}$ generated by $h$ 
is a positive-definite lattice of rank $1$ by $b_X$.
Let $p_1, \dots, p_t$ be the singular points of $Y$.
For each $p_i$, 
we denote by $\Sigma_i\subset H_2(X)$
the submodule generated by the homology classes of
integral curves
on $X$ that are contracted to $p_i$ by $\rho$.
By the theorem of Mumford~\cite{MR0153682},
the $\Z$-module $\Sigma_i$  is a negative-definite lattice by $b_X$.
The lattice $\Sigma_i$ is perpendicular to $\ang{h}$ and $\Sigma_j$ $(j\ne i)$
with respect to $b_X$.
Therefore the submodule 
$$
\MC\sp 0:=\ang{h}\perp\Sigma_1\perp\dots\perp \Sigma_t
$$
of $H_2(X)$ is a lattice by $b_X$.
Let $C_i$ be an irreducible component of $C$,
and let $\wt{C}_i$ be the integral curve on $X$
such that $\phi (\wt{C}_i)=C_i$.
Since $\Pic(S)\cong \Z$
is generated by the class of $\HHH$,
there exists an integer $d_i$ such that 
$C_i$ is linearly equivalent to
$\HHH\sp{\otimes d_i}$ on $S$.
Then 
there exists $\gamma \in \Sigma_1\perp\dots\perp \Sigma_t$
such that 
$2[\wt{C}_i]=d_i h +\gamma$
holds in $H_2(X)$.
Therefore $\MCspsim \otimes \Q$ is equal to $\MC\sp 0\otimes \Q$ in $H_2(X)\otimes \Q$.
\end{proof}
From now on, we consider $N_C$ as a lattice by $b_X$.
\begin{lemma}\label{lem:choice}
The isomorphism class of the lattice $\NC$ does not depend on the choice 
of the morphism  $\rho : X\to Y$.
\end{lemma}
\begin{proof}
Let $X\sprime\to X$ be the blowing up  
at a point $P$ on $\phi\inv (C)$,
and let $E$ be the $(-1)$-curve on $X\sprime$ contracted to $P$.
Then we have a natural isomorphism 
$$
H_2 (X\sprime) = H_2 (X)\perp\ang{[E]}.
$$
Hence the lattice $\NC\sprime\subset H_2(X\sprime)\tf$ constructed from  $X\sprime$
by the method described above  is isomorphic to $\NC$.
\end{proof}
Therefore we can consider the isomorphism class of the lattice $\NC$ as an invariant of
the $\LLL$-curve $C$.
\begin{definition}\label{def:Gamma}
Let $C_1, \dots, C_m$ be the irreducible components of an $\LLL$-curve $C$.
We denote by $\Gamma_i\subset \pione (S\sm C)$ the conjugacy class of 
simple loops around $C_i$, and put
$\varGamma (C):=\{\Gamma_1, \dots, \Gamma_m\}$.
Let $C$ and $C\sprime$ be $\LLL$-curves.
We say that $(S, C)$ and $(S, C\sprime)$   are  \emph{$\varGamma$-equivalent}
if there exists a homeomorphism $S\sm C\cong S\sm C\sprime$
such that the induced 
 isomorphism $\pione (S\sm C)\cong \pione (S\sm C\sprime)$
 gives rise to  a bijection from $\varGamma (C)$ to $\varGamma (C\sprime)$.
\end{definition}
It is obvious that, for  $\LLL$-curves $C$ and $C\sprime$,
if $(S, C)$ and $(S, C\sprime)$ are homeomorphic,
then  they  are $\varGamma$-equivalent.
\par
\medskip
The following is the main result of this paper:
\begin{theorem}\label{thm:main}
The isomorphism class of the lattice $\NC$ is an invariant of the $\varGamma$-equivalence class of $(S, C)$.
\end{theorem}
\begin{proof}
Since $\pione (S)$ is assumed to be trivial,
$\pione (S\sm C)$ is generated by  simple loops
around  irreducible components of $C$.
Therefore 
a homomorphism
$\pione (S\sm C)\to\Z/2\Z$
that maps every element of $\Gamma_1\,\cup\dots\cup \,\Gamma_m$
to the non-trivial element of $\Z/2\Z$
is unique.
Consequently,   the homeomorphism type of
$$
U:=\phi\inv (S\sm C)\;\;\subset\;\; X
$$
is uniquely determined by the   $\varGamma$-equivalence class of $(S, C)$.
For a compact subset $K$ of $U$,
we denote by $B^K$ the image of the natural homomorphism
$H_2(U\sm K)\to H_2(U)$.
We then put
$$
B_\infty \;\;:=\;\;\bigcap B^K,
$$
where $K$ runs through the set of all compact subsets of $U$,
and set 
$$
\BUspsim:=H_2(U)/B_\infty
\quand
\BU:=(\BUspsim)\tf.
$$
Since every topological cycle is compact, 
the intersection pairing
$b_U$ on $H_2(U)$ sets up  a symmetric bilinear form
$$
\map{\beta_U}{\BUspsim\times\BUspsim}{\Z}.
$$
By construction,
the isomorphism class of $(\BU, \beta_U)$
is determined by the homeomorphism type of $U$,
and hence by the  $\varGamma$-equivalence class of $(S, C)$.
Therefore  it is enough to show that 
the lattice $\NC$ is isomorphic to $(\BU, \beta_U)$.
\par
\medskip
We put $D:=\phi\inv (C)$,
and equip $D$ with the reduced structure.
Let $T\subset X$ be a tubular neighborhood of  $D$.
We  put
$$
T\sptimes :=T\sm D=T\cap U,
$$
and consider the Mayer-Vietoris sequence
 $$
\renewcommand{\arraystretch}{1.2}
\begin{array}{ccccccc}
\maprightsp{} & H_2 (T\sptimes) & \maprightsp{i} & H_2(T)\oplus H_2(U) &\maprightsp{j} & H_2(X) & \maprightsp{}\mystrutd{10pt}\\
& x& \longmapsto & (i_T(x), i_U(x)) &&&\\
&&&(y, z)&\longmapsto & j_T(y)-j_U(z) & .\\
\end{array}
$$
First we prove
\begin{equation}\label{eq:ImjU}
\Im (j_U)=\NCspsim .
\end{equation}
It is obvious that $\Im (j_U)\subseteq \NCspsim $.
Let $[W]$ be an element of $\NCspsim $
represented by a $2$-dimensional topological cycle $W\subset X$.
We can assume the following:
\begin{rmenumerate}
\item $W\cap \Sing (D)=\emptyset$, 
\item $W\cap D$ consists of a finite number of points, and
\item locally around each point $P$ of $W\cap D$,
$W$ is a $\Cinf$-manifold intersecting $D$ at $P$ transversely.
\end{rmenumerate}
Let $D^{(1)}, \dots, D^{(n)}$ be the irreducible   components of $D$.
For each $\nu=1, \dots, n$,
let $P_{+, 1}^{(\nu)}, \dots, P_{+, k(\nu)}^{(\nu)}$ 
(resp.~$P_{-, 1}^{(\nu)}, \dots, P_{-, l(\nu)}^{(\nu)}$)
be the intersection points of $W$ and $D^{(\nu)}$ with 
local intersection number $1$ (resp.~$-1$).
Since $b_X([W], [D\spnu])=0$,
we have 
$$
k(\nu)=l(\nu).
$$
Let $\Delta\subset \C$ be the closed unit disk,
and let $I\subset\R$ be the closed interval $[0, 1]$.
Since $D\spnu\sm\Sing (D)$ is path-connected for each $\nu$, 
we have  continuous maps
$$
\map{\xi_i\spnu}{\Delta\times I}{X}
$$
for $\nu=1, \dots, n$ and $i=1, \dots, k(\nu)$
with the following properties:
\begin{rmenumerate}
\item $\xi_i\spnu (\Delta\times I) \cap D \;\;\subset\;\; D\spnu\sm\Sing (D)$, 
\item $(\xi_i\spnu)\inv (D)=\{0\}\times I$,
\item $\xi_i\spnu(0,0)=P_{+, i}\spnu$, and $\xi_i\spnu$ induces a homeomorphism from $\Delta\times\{0\}$ to
a closed neighborhood $\Delta_{+, i}\spnu\subset W$ of $P_{+, i}\spnu$ in $W$, and 
\item $\xi_i\spnu(0,1)=P_{-, i}\spnu$, and $\xi_i\spnu$ induces a homeomorphism from $\Delta\times\{1\}$ to
a closed neighborhood $\Delta_{-, i}\spnu \subset W$ of $P_{-, i}\spnu$ in $W$.
\end{rmenumerate}
We then put
$$
W\sprime:=\Bigl( W\sm \bigcup_{\nu, i} \,(\Delta_{+, i}\spnu\cup \Delta_{-, i}\spnu) \Bigr)\; \cup \;
\bigcup_{\nu, i} \, \xi_i\spnu (\partial\Delta\times I).
$$
That is, we cut out discs  $\Delta_{+, i}\spnu$ and  $\Delta_{-, i}\spnu$ from  $W$,
and put tubes
$\xi_i\spnu (\partial\Delta\times I)$
that connect pairs of circles 
$\partial \Delta_{+, i}\spnu$ and $\partial  \Delta_{-, i}\spnu$
around  the paths $\xi_i\spnu(\{0\}\times I)$ on $D\spnu$
from $P\spnu_{+, i}$ to $P\spnu_{-, i}$.
Since the local intersection numbers of $W$ and $D\spnu$
at $P\spnu_{+, i}$ and at $P\spnu_{-, i}$
have opposite signs,
we can put an orientation on each   $\xi_i\spnu (\Delta\times I)$
in such a way that
$\Delta\spnu_{+, i}\subset W$ and
$\xi_i\spnu (\Delta\times \{0\})\subset \partial (\xi_i\spnu (\Delta\times I))$
(resp. 
$\Delta\spnu_{-, i}\subset W$ and 
$\xi_i\spnu (\Delta\times \{1\})\subset \partial (\xi_i\spnu (\Delta\times I))$
)
have the opposite orientations.
Then, with the orientation on 
the tubes $\xi_i\spnu (\partial\Delta\times I)$
induced from the orientation of $\xi_i\spnu (\Delta\times I)$,
the space 
 $W\sprime$  becomes a topological cycle.
Note that  $W$ and $W\sprime$ are
homologous in $X$, because $W- W\sprime$
is the boundary of the $3$-dimensional topological chain $\bigcup \, \xi_i\spnu (\Delta\times I)$.
Moreover $W\sprime$ is disjoint from $D$.
Therefore $[W]=[W\sprime]$ is contained in $\Im (j_U)$, and hence~\eqref{eq:ImjU} is proved.
\par
\medskip
Since $T$ is a  tubular neighborhood of $D$, we have
\begin{equation}\label{eq:ImiU}
\Im (i_U)=B_\infty. 
\end{equation}
If $j_U(z)=0$, then
$(0, z)\in \Ker (j)=\Im (i)$,
and hence $z\in \Im (i_U)$.
Therefore we have
$\Ker (j_U)\subseteq \Im (i_U)$.
Consider the following commutative diagram:
\newcommand{\downsurj}{\raisebox{10pt}{\rotatebox{-90}{$\surj$}}}@
\newcommand{\downinj}{\raisebox{8pt}{\rotatebox{-90}{$\inj$}}}
\newcommand{\verteq}{\raisebox{-2pt}{\rotatebox{90}{$=$}}}@
\begin{equation}\label{eq:diag}
\renewcommand{\arraystretch}{1.4}
\begin{array}{ccccccccc}
0 & \maprightsp{} & \Ker (j_U) & \maprightsp{} &H_2 (U) & \maprightsp{j_U} & \NCspsim  & \maprightsp{} & 0\phantom{.}\\
 && \downinj & &\verteq & & \;\;\,\downsurj\;\;\raise 3pt \hbox{${}_\lambda$} &&\\
 0 & \maprightsp{} & \Im (i_U) &\maprightsp{} & H_2 (U) & \maprightsp{} & \BUspsim  & \maprightsp{} & 0.
\end{array}
\end{equation}
The upper sequence is exact by~\eqref{eq:ImjU},
and the lower sequence is exact by~\eqref{eq:ImiU}.
Remark that, by the definition of the intersection pairing,
we have
\begin{equation*}
b_U(z, z\sprime)=b_X(j_U(z), j_U(z\sprime))
\end{equation*}
for any $z, z\sprime \in H_2(U)$.
Hence, for any $\zeta, \zeta\sprime\in \NCspsim\subset H_2(X)$, we have
\begin{equation}\label{eq:bbl}
b_X(\zeta, \zeta\sprime)=\beta_U (\lambda(\zeta), \lambda(\zeta\sprime)),
\end{equation}
where $\lambda: \NCspsim\to \BUspsim$ is the vertical surjection in the diagram~\eqref{eq:diag}.
Therefore,
in order to show that
$N_C$ is isomorphic to $(\BU, \beta_U)$,
 it is enough to prove that $\lambda$ induces
an isomorphism $\NC\isom \BU$ on the torsion-free quotients,
or equivalently, $\lambda\inv ( (\BUspsim)\tor)=(\NCspsim)\tor$ holds.
It is obvious that $\lambda((\NCspsim)\tor)\subseteq (\BUspsim)\tor$.
Suppose that $\zeta\in \NCspsim$ satisfies $\lambda(\zeta)\in (\BUspsim)\tor$.
By~\eqref{eq:bbl},
we have
$b_X(\zeta, \zeta\sprime)=0$ for any $\zeta\sprime\in \NCspsim$.
Since $\NC$ is a lattice,
we have $\zeta\in (\NCspsim)\tor$.
\end{proof}
\begin{corollary}\label{cor:PC}
Let $C$ and $C\sprime$ be plane curves of the same degree.
Suppose that $\deg C=\deg C\sprime$ is even.
If $(\Pt, C)$ and $(\Pt, C\sprime)$ are homeomorphic,
then $N_C$ and $N_{C\sprime}$ are isomorphic.
\end{corollary}
\section{Sextics with only simple singularities}\label{sec:conn}
Let $\P_*(H^0(\Pt, \OOO_{\Pt}(6)))$
be the projective space
of one-dimensional subspaces
of the vector space $H^0(\Pt, \OOO_{\Pt}(6))$
of homogeneous polynomials of degree $6$ on $\Pt$.
For a Dynkin type $R$ of rank $\le 19$,
we denote by 
$$
\MMM (R)\;\;\subset\;\;\P_*(H^0(\Pt, \OOO_{\Pt}(6)))
$$
the space of
$ADE$-sextics   of type $R$.
Using Urabe's idea~\cite{MR1101859}, 
Yang~\cite{MR1387816} made the complete list of  Dynkin types $R$
such that $\MMM (R)\ne \emptyset$.
Degtyarev~\cite{degtyarev-2005} refined Yang's argument,
and gave a method to calculate the connected components of  $\MMM(R)$
for a given $R$.
In this section,
we expound Degtyarev's theory.
\par
\medskip
We fix some notation and terminologies about lattices.
\par
\medskip
Let $\lat$ be a lattice of rank $n=2 + s_-$ and 
signature $(2, s_-)$.
For a non-zero vector $\omega \in \lat\tC$,
we denote by $[\omega]\in \P_* (\lat\tC)$
the one-dimensional vector space spanned by $\omega$.
We put
$$
\Omega_\lat:=\set{[\omega]\in \P_*(\lat\tC)}{(\omega, \omega)=0, \,(\omega, \bar\omega)>0}.
$$
It is easy to verify that  $\Omega_\lat$ is  a complex manifold of dimension $s_-=n-2$
consisting of     two connected components.
\par\smallskip
The \emph{dual lattice $\lat\dual$} of a lattice $\lat$ is defined by
$$
\lat\dual:=\set{v\in \lat\tQ}{(x, v)\in \Z\;\;\textrm{for all}\;\; x\in \lat}.
$$
We have $\lat\subset \lat\dual$.
An \emph{overlattice} of $\lat$ is a submodule $\lat\sprime$ of $\lat\dual$ containing $\lat$
such that the natural $\Q$-valued symmetric bilinear form on $\lat\tQ$
  takes values in 
$\Z$ on $\lat\sprime$.
The \emph{discriminant group} $G_\lat$ of $\lat$ is defined by 
$$
G_\lat:=\lat\dual/\lat.
$$
A lattice is called \emph{unimodular} if $\lat\dual=\lat$.
A lattice $\lat$ is said to be \emph{even} if $(v, v)\in 2\,\Z$ 
holds for every $v\in \lat$.
If $\lat$ is an even lattice,   we can define a quadratic form
$$
q_\lat\;:\; G_\lat\to\Q/2\,\Z
$$
by $q_\lat (\bar v):=(v, v)\,\bmod 2\,\Z$,
where $v\in \lat\dual$ and $\bar v:= v\,\bmod \lat$.
This quadratic form is called the \emph{discriminant form} of  $\lat$.
See Nikulin~\cite{MR525944} for the  basic properties of discriminant forms.
\par
\smallskip
Let $\lat$ be a negative-definite even lattice.
A vector $d\in \lat$ is called a \emph{root} if $(d, d)=-2$ holds.
We say that $\lat$ is a \emph{root lattice} if $\lat$ is generated by the roots in $\lat$.
The isomorphism classes of root lattices  are  in one-to-one correspondence
with the Dynkin types
(see, for example, Ebeling~\cite[Section 1.4]{MR1938666}).
We denote by  $\rootlat_R^-$
 the negative-definite root lattice of Dynkin type $R$.
A subset $F$ of $\rootlat_R^-$
is called a \emph{fundamental system of roots}
if every element of $F$ is a root, $F$ is a basis of $\rootlat_R^-$,  and 
each  root $v\in \rootlat_R^-$ is  written as a linear combination 
$v=\sum_{d\in F}  k_d d$ 
of roots $d\in F$ with integer coefficients $k_d$ all non-positive or all non-negative.
A fundamental system $F$ of roots
 exists
 (see~\cite[Section 1.4]{MR1938666}).
 The intersection matrix of $\rootlat_R^-$  with respect to the basis $F$
is the Cartan matrix of type $R$ multiplied by $-1$.
\par
\smallskip
A lattice is called a \emph{$K3$ lattice} if it is even, unimodular,
of rank $22$ and  with signature $(3, 19)$.
By the structure theorem of unimodular lattices,
a $K3$ lattice is unique up to isomorphisms
(see,  for example, Serre~\cite[Chapter V]{MR0344216}).
\par
\medskip
We now start explaining Degtyarev's theory.
Let $R$ be a Dynkin type with
$$
r:=\rank (R)\le 19.
$$
\par
\smallskip
First, we define a set $Q(R)$ and an equivalence relation $\relQ$ on it.
We denote by 
$\ang{h}$ the lattice of rank $1$ generated by a vector $h$ with $(h, h)=2$.
We put
$$
M^0:=\rootlat_R^{-}\perp  \ang{h},
$$
which is an even lattice of signature $(1, r)$.
\emph{We choose a fundamental system of roots $F\subset \rootlat_R^{-}$
once and for all},
and put
$$
\AutFh:=\set{g\in O(M^0)}{g(F)=F, \;g(h)=h}.
$$
We denote by $\Ms$ the set of even overlattices $M$ of $M^0$
satisfying the following two conditions:
\begin{itemize}
\item[(m1)]
$\set{v\in M}{(v, h)=1, \;(v, v)=0}=\emptyset$, and 
\item[(m2)]
$\set{v\in M}{(v, h)=0, \;(v, v)=-2}=\set{v\in \rootlat_R^{-}}{(v, v)=-2}$.
\end{itemize}
(These conditions correspond to the conditions (a) and (b) in~\cite[Theorem 2.3]{MR1387816}.)
For $M\in \Ms$, we denote by $\Ns (M)$
a complete set of representatives of isomorphism classes of even lattices $N$
of rank $21-r$
satisfying the following two conditions:
\begin{itemize}
\item[(n1)]
$N$ is  of signature $(2, 19-r)$, and
\item[(n2)]
the discriminant form 
$(G_N, q_N)$ of $N$ is isomorphic to $(G_M, -q_M)$.
\end{itemize}
By Nikulin~\cite[Proposition~1.6.1]{MR525944},
the conditions (n1) and  (n2) are  equivalent to the following condition:
\begin{itemize}
\item[(n)]there exists an even unimodular overlattice $L$ of $M\perp N$
with signature $(3, 19)$ such that 
  $M$ and $N$ are primitive  in $L$.
\end{itemize}
Let $N$ be an element of $\Ns(M)$.
We denote by $\Ls(M, N)$
the set of even unimodular overlattices $L$ of $M\perp N$
such that $M$ and $N$ are primitive  in $L$.
Note that 
every  $L\in \Ls (M, N)$ is a $K3$ lattice.
We also denote by $\cpersps (N)$ the set of connected components  
of the complex manifold $\Omega_N$.
Remark that  we have   $|\cpersps (N)|=2$.
\par
\smallskip
We define $Q(R)$ to be the set of quartets $(M, N, L, \cpersp)$
such that $M\in \Ms$, $N\in \Ns(M)$, $L\in \Ls (M, N)$,
and $\cpersp\in \cpersps (N)$.
For quartets $(M, N, L, \cpersp)$ and 
$(M\sprime, N\sprime, L\sprime, \cpersp\sprime)$ in $Q(R)$,
we write
$$
(M, N, L, \cpersp)\;\;\relQ \;\;(M\sprime, N\sprime, L\sprime, \cpersp\sprime)
$$
if the following hold.
\begin{itemize}
\item[(i)]
There exists $g^0\in \AutFh\subset O(M^0)$ such that 
the induced action of $g^0$ on $\Ms$ maps $M\in\Ms$
to $M\sprime\in \Ms$.
We denote by
$g_M: M\isom M\sprime$
the unique isomorphism  satisfying  $g_M|M^0=g^0$.
\item[(ii)]
Since $(G_M, -q_M)$ and $(G_{M\sprime}, -q_{M\sprime})$ are isomorphic,
there exists a canonical bijection between 
$\Ns (M)$ and $\Ns(M\sprime)$.
The elements $N\in \Ns(M)$ and $N\sprime\in \Ns(M\sprime)$ are 
corresponding by this bijection; that is,  $N$ and $N\sprime$ are isomorphic.
\item[(iii)]
There exists an isomorphism 
$g_N: N\isom N\sprime$
of lattices 
such that 
the bijection 
$\Ls (M, N)\isom \Ls (M\sprime, N\sprime)$
induced by the isomorphism $g_M\oplus g_N$
from $ M\perp N $ to $ M\sprime\perp N\sprime$
 maps $L\in \Ls (M, N)$ to $L\sprime\in \Ls (M\sprime, N\sprime)$,
and that the induced isomorphism $\Omega_{N}\isom \Omega_{N\sprime}$
maps $\cpersp$ to $\cpersp\sprime$.
\end{itemize}
For $(M, N, L, \cpersp)\in Q(R)$,
we denote by $[M, N, L, \cpersp]\in Q(R)/\relQ$
the equivalence class of the relation $\relQ$ containing $(M, N, L, \cpersp)$.
\begin{remark}\label{rem:rho1}
If $(M,N,L,\cpersp)\in Q(R)$,
then 
$M^0$ is the sublattice of $L$
generated by $F\subset L$ and $h\in L$,
$M$ is the primitive closure of $M^0$ in  $L$,
and $N$ is the orthogonal complement of $M$ in $L$.
Hence 
$[M,N,L, \cpersp]=[M\sprime,N\sprime,L\sprime,\cpersp\sprime]$
holds if and only if 
there exists an isomorphism $L\isom  L\sprime$
that maps $F$ to $F$, $h$ to $h$,
and 
such that
the induced isomorphism  $\Omega_L\isom \Omega_{L\sprime}$
maps 
the connected component $\cpersp$ of $\Omega_N\subset \Omega_L$
to the connected component $\cpersp\sprime$
of $\Omega_{N\sprime}\subset \Omega_{L\sprime}$.
\end{remark}
Next
we define a map $\rho$ from the space $\MMM(R)$
to the set  $Q(R)/\relQ$.
Let $C$ be an $ADE$-sextic of type $R$,
and let $X$ be the minimal resolution of the double covering $Y\to\Pt$
that branches exactly  along $C$.
We denote by $\polar$ the line bundle on $X$
corresponding to the pull-back of the invertible sheaf $\OOO_{\Pt} (1)$.
We have $([\polar], [\polar])=2$. We then put
$$
L_X:=H^2(X, \Z),
$$
which is a $K3$ lattice.
Let $F\sbXL  \subset L_X$ be 
the set of cohomology classes of $(-2)$-curves that are contracted by 
the desingularization  morphism  $X\to Y$,
and  let $\rootlat\sbXL\subset L_X$ be 
 the sublattice of $L_X$ generated by $F\sbXL  $.
Then 
$\rootlat\sbXL$ is a negative-definite root lattice of type $R$.
It is known that 
 $F\sbXL  $ is a fundamental system of roots in $\rootlat\sbXL$
 (see~\cite[Proposition 2.4]{normalK3}).
 In particular,
 there exists an isomorphism of lattices from $\rootlat\sbXL$ to $\rootlat^-_R$
 that maps  $F\sbXL  $  to the fixed fundamental system of roots $F\subset \rootlat^-_R$ bijectively.
We put
$$
M^0\sbXL:=\rootlat\sbXL\perp\ang{[\polar]}, 
$$
and  choose an isomorphism
$$
\mapisom{\gamma^0_M}{M^0\sbXL}{M^0}
$$
satisfying $\gamma^0_M (F\sbXL  )=F$ and $\gamma^0_M([\polar])=h$.
Let $M\sbXL$ be the primitive closure of $M^0\sbXL$
in $L_X$, and 
$M$  the  even overlattice of $M^0$
corresponding to 
the even overlattice 
$M\sbXL$ of $M^0\sbXL$ by $\gamma^0_M$.
Then $M$ satisfies the conditions  (m1)  and  (m2)
(see~\cite[Proposition 2.1]{normalK3}).
Hence $M\in \Ms$.
We denote by 
$$
\mapisom{\gamma_M}{ M\sbXL}{ M}
$$
the isomorphism induced by $\gamma^0_M$.
Let $N\sbXL$ be the orthogonal complement of $M\sbXL$ in $L_X$.
Since the $K3$ lattice $L_X$ is an even unimodular overlattice of $M\sbXL\perp N\sbXL$
in which $M\sbXL$ and $N\sbXL$ are  primitive,
the lattice $N\sbXL$ satisfies the condition (n).
Hence there exists a unique element $N$ of $\Ns (M)$ that is isomorphic to $N\sbXL$.
We choose an isomorphism
$$
\mapisom{\gamma_N }{ N\sbXL}{ N}.
$$
By the isomorphism
$$
\mapisom{\gamma_M\oplus  \gamma_N\,}{\, M\sbXL\perp N\sbXL \,}{\, M\perp N},
$$
the  even unimodular overlattice $L_X$ of $M\sbXL\perp N\sbXL$
corresponds 
to an element $L$ of $\Ls (M, N)$.
We denote by
$$
\omega_X\in H^{2, 0}(X)\subset L_X\tC
$$
the cohomology class of 
a non-zero holomorphic $2$-form on $X$.
Since $M\sbXL\subset H^{1,1}(X)$,
the vector $\omega_X$ defines a point $[\omega_X]$ of $\Omega_{N\sbXL}$.
Let $\cpersp$ be the connected component of $\Omega_N$
that contains the point $[\gamma_N(\omega_X)]$.
Thus we obtain a quartet $(M, N, L, \cpersp)\in Q(R)$.
The choices we have made during the process 
of finding $(M, N, L, \cpersp)$ are only on
$\gamma^0_M$ and $\gamma_N$.
Since $\gamma^0_M$ is unique up to $\AutFh$ and
$\gamma_N$ is unique up to $O(N)$, 
 the equivalence class $[M, N, L, \cpersp]$
does not depend on these choices.
We thus  can put
$$
\rho (C) :=[M, N, L, \cpersp].
$$
\begin{remark}\label{rem:rho2}
By definition,
we have $\rho (C)=[M, N, L, \cpersp]$
if and only if there exists an isomorphism
$L_X\isom L$
that maps $F\sbXL  $ to  $F$,
$[\polar]$ to $h$,
and such that the induced isomorphism $\Omega_{L_X}\isom\Omega_L$
maps
the point $[\omega_X]\in \Omega_{N\sbXL} \st\Omega_{L_X}$ to a point of 
the connected component $\cpersp$ of $\Omega_N\subset \Omega_L$.
\end{remark}
We now have all the ingredients that are needed to state the main theorem 
of Degtyarev~\cite{degtyarev-2005}:
\begin{theorem}
The map $\rho$ induces a bijection from 
the set of connected components of the space $\MMM(R)$  to the set $Q(R)/\relQ$.
\end{theorem}
The main tool of the proof is
the Torelli theorem for the refined period map
of marked $K3$ surfaces.
See the book by Barth, Hulek, Peters and Van~de Ven~\cite[Theorems 12.3 and 14.1 in Chapter VIII ]{MR2030225}.
\par
\medskip
By definition, we have the following:
\begin{corollary}\label{cor:main}
Let $C$ be an $ADE$-sextic such  that $\rho (C)=[M,N,L,\cpersp]$.
Then the lattice $N$ is isomorphic 
to  the invariant $N_C$ of the $\varGamma$-equivalence class of $(\Pt, C)$. 
\end{corollary}
We explain how to calculate the set $Q(R)/\relQ$.
By~\cite[Proposition 1.4.1]{MR525944},
the  even overlattices of $M^0=\rootlat_{R}^-\perp \ang{h}$
are in one-to-one correspondence
with the  totally isotropic subgroups
of the discriminant form $(G_{M^0}, q_{M^0})$.
For an even  overlattice $M$ of $M^0$,
we can determine  whether $M$
satisfies the conditions (m1) and (m2)
by the method described in~\cite{MR2036331}.
Since $G_{M^0}$ is finite,
we obtain the set $\Ms$.
The group $\AutFh$
is isomorphic to the automorphism group
of the Dynkin diagram of type $R$,
and hence it is finite.
Therefore 
the image of the natural homomorphism 
$$
\AutFh\;\;\inj \;\; O(M^0)\;\;\to \;\; O(q_{M^0})
$$
is  easy to calculate,
where $O(q_{M^0})$ is the automorphism group
of the finite quadratic form $(G_{M^0}, q_{M^0})$
(see~\cite[Section 6.2]{MR1813537}).
Consequently
 we obtain the set
$$
\overline{\Ms}:=\AutFh\backslash \Ms
$$
of the orbits of the action of $\AutFh$ on $\Ms$.
For an element $M$ of $\Ms$,
let $[M]\in \overline{\Ms}$
denote the orbit containing $M$.
We  also put
$$
\AutFhM :=\set{g\in \AutFh}{\textrm{$g$ fixes $M\in \Ms$}}.
$$
We  have a natural map
$$
\pr \;:\; Q(R) /\relQ\;\;\to \;\;\overline{\Ms}
$$
that maps $[M,N,L,\cpersp]$ to $[M]$.
We denote by $\overline{\Ms}\spsharp\subset\overline{\Ms}$ the image of the map $\pr : Q(R) /\relQ\to \overline{\Ms}$\,;
that is, we put
$$
\overline{\Ms}\spsharp:=\set{[M]\in\overline{\Ms}}{\Ns(M)\ne\emptyset}.
$$
For $[M]\in \overline{\Ms}$, we can determine whether $\Ns (M)$ is empty or not by the criterion of 
Nikulin~\cite[Theorem 1.10.1]{MR525944}.
Hence $\overline{\Ms}\spsharp$ is calculated.
\par
\medskip
Suppose that $[M]\in \overline{\Ms}\spsharp$.
By~\cite[Corollary 1.9.4]{MR525944}, 
the set $\Ns (M)$ forms  a genus.
If $r:=\rank (R)<19$,
then the isomorphism class of an indefinite  lattice $N$ of signature $(2, 19-r)$ 
is determined by the spinor genus by Eichler's theorem
(see, for example, Cassels~\cite{MR522835}).
The method of  enumeration  of spinor genera in a given genus
is  described in Conway and Sloane~\cite[Chapter 15]{MR1662447}.
When  $\rank (R)=19$,
the set $\Ns(M)$ is easily calculated by Corollary~\ref{cor:BF} below.
\par
\medskip
For each $[M]\in \overline{\Ms}\spsharp$, 
we have a natural map
$$
\pr_{[M]} \;:\; \pr\inv([M])\;\;\to \;\;\Ns(M)
$$
that maps $[M\sprime,N\sprime,L\sprime,\cpersp\sprime]\in \pr\inv([M])$ to 
the lattice $N\in \Ns(M)$ isomorphic  to $N\sprime\in \Ns(M\sprime)$.
(Note that, if $[M]=[M\sprime]$,
then $M$ and $M\sprime$ are isomorphic, 
and hence $\Ns(M)$ and $\Ns(M\sprime)$ are canonically identified.)
Let $N$ be an element of $\Ns(M)$.
We put
$$
F([M], N):=\pr_{[M]}\inv (N).
$$
We can regard $\AutFhM $ as a subgroup of $O(M)$:
$$
\AutFhM=\set{g\in O(M)}{g(F)=F, \, g(h)=h}.
$$
Then the group
$\AutFhM \times O(N)$
acts on the set $\Ls (M, N)\times \cpersps (N)$
in the natural way.
The fiber $F([M], N)$ of $\pr_{[M]}$ over $N$ is, by definition,
equal to the set of orbits of this action:
$$
F([M], N)=(\AutFhM \times O(N))\backslash (\Ls (M, N)\times \cpersps (N)).
$$
By~\cite[Proposition~1.6.1]{MR525944}, 
there exists a natural bijection between 
 the set $\Ls(M, N)$ 
 and the set of isomorphisms of finite quadratic forms 
from 
$(G_M, -q_M)$ to $(G_N, q_N)$.
Since $G_M\cong G_N$ is a finite abelian group,
 we obtain the set $\Ls(M, N)$.
Hence the set $F([M], N)$ can be  calculated,
provided that the group $O(N)$ and its actions on $(G_N, q_N)$ and on $\cpersps(N)$
are calculated.
\begin{remark}
When $\rank (R)=19$,
the lattice $N$ is  positive-definite. 
Hence $O(N)$ is finite,
and we can easily make the list of elements of $O(N)$. 
The  actions of $O(N)$ on $(G_N, q_N)$ and on $\cpersps(N)$
are  then readily   calculated.
\end{remark}
We use the following terminology in \S\ref{sec:examples} and \S\ref{sec:singK3}.
\begin{definition}
Let  $\tau : \cpersps(N)\isom  \cpersps(N)$
be the transposition  of the two connected components of $\Omega_N$.
An orbit $\Orb\subset \Ls(M, N)\times\cpersps(N)$
of the action of $\AutFhM \times O(N)$
is called \emph{real}
if $\Orb$ is stable under $\tau$.
\end{definition}
We review  the classical theory of
binary forms due to Gauss
(see Edwards~\cite{MR1416327} or 
Conway and Sloane~\cite[Chapter 15]{MR1662447}).
For integers $a, b, c$,
we denote by
$\Mat [a,b,c]$ 
the  matrix
$$
\left[\begin{array}{cc} a & b \\ b & c \end{array}\right].
$$
For a positive integer $d$, we put
$$
\Mats_d:=\set{\Mat[a,b,c]}{a\equiv c\equiv  0 \,(\bmod\,  2),\;\; a>0, \;\; c>0, \;\; ac-b^2=d},
$$
on which $\GL_2(\Z)$ acts from right by $(Q, g)\mapsto {}^T \hskip -1pt g \,Q\, g$.
The set of isomorphism classes of even positive-definite lattices of rank $2$ with discriminant
$d$ is canonically identified with
the set $\Mats_d/\mathord{\GL}_2(\Z)$  of $\GL_2(\Z)$-orbits in $\Mats_d$.
\begin{definition}
We call an $\SL_2(\Z)$-orbit in $\Mats_d$
an isomorphism class of even positive-definite \emph{oriented} lattices 
of rank $2$ with discriminant $d$.
\end{definition}
For $\Mat[a,b,c]\in\Mats_d$,
we denote by $\Lat[a,b,c]$ (resp.~$\oriLat[a,b,c]$)
the lattice (resp.~the oriented lattice)
expressed  by $Q[a,b,c]$.
\begin{proposition}\label{prop:BF}
Let $d$ be a positive integer.
Then the  set
$$
\set{\oriLat  [a,b,c]}{\Mat[a,b,c]\in\Mats_d,\;\; 
-a<2b\le a\le c \;\;{\rm with}\;\; b\ge 0 \;\;{\rm if}\;\;a=c}
$$
is a complete set  of representatives of isomorphism classes of 
even positive-definite oriented  lattices of rank $2$
with discriminant $d$.
\end{proposition}
\begin{corollary}\label{cor:BF}
Let $d$ be a positive integer.
Then the  set
\begin{equation}\label{eq:QGL}
\set{\Lat  [a,b,c]}{\Mat[a,b,c]\in\Mats_d,\; 0\le 2b\le a\le c}
\end{equation}
is a complete set  of representatives of isomorphism classes of 
even positive-definite   lattices of rank $2$
with discriminant $d$.
\end{corollary}
\begin{remark}\label{rem:forget}
Let $\Lat[a,b,c]$ be an element of the set~\eqref{eq:QGL}, and 
let $[\Lat[a,b,c]]\in  \Mats_d/\mathord{\GL}_2(\Z)$ be 
the  $\GL_2(\Z)$-orbit containing $\Lat[a,b,c]$.
Then the fiber of 
the natural map
$\Mats_d/\mathord{\SL}_2(\Z)\to \Mats_d/\mathord{\GL}_2(\Z)$
over $[\Lat[a,b,c]]$
consists of two  elements if 
$$
0< 2b< a< c,
$$ 
while it consists of a single element  otherwise.
\end{remark}
\section{Examples of arithmetic Zariski pairs}\label{sec:examples}
Let $f\in \Q[t]$ be 
an irreducible polynomial.
We denote by $F_f$ the field $\Q[t]/(f)$.
Suppose that a homogeneous polynomial $\Phi(x_0,x_1,x_2)$
with coefficients in $F_f$ is given.
For a complex root $\alpha$ of $f$,
we denote by $\sigma_\alpha: F_f\inj \C$
the embedding given by $t\mapsto \alpha$,
and by $\Phi^\alpha$ the homogeneous polynomial
obtained by applying $\sigma_\alpha$ to the coefficients of $\Phi$.
We say that two plane curves $C$ and $C\sprime$ are \emph{conjugate in $F_f$}
if there exist a homogeneous polynomial $\Phi$
with coefficients in $F_f$ and distinct complex roots $\alpha$ and $\alpha\sprime$
of $f$ such that
$C=\{\Phi^\alpha=0\}$ and $C\sprime=\{\Phi^{\alpha\sprime}=0\}$.
\par
\medskip
Suppose that $C$ and $C\sprime$ are conjugate $ADE$-sextics.
Then the configurations of $C$ and $C\sprime$
are the same.
(See Yang~\cite[\S3]{MR1387816}
for the precise definition of the configuration
of an $ADE$-sextic.)
In particular, there exist  tubular neighborhoods $T\subset \Pt$ of $C$ and $T\sprime\subset \Pt$ of $C\sprime$
such that $(T, C)$ and $(T\sprime, C\sprime)$ are diffeomorphic.
\par
\medskip
Combining this fact with Corollaries~\ref{cor:PC} and~\ref{cor:main},
we see that the following pairs of conjugate maximizing sextics discovered by
Artal-Bartolo, Carmona-Ruber and  Cogolludo-Agust{\'{\i}}n~\cite[Theorem~5.8]{MR1900779}
are in fact arithmetic Zariski pairs.
\begin{example}\label{example:1}
Consider the Dynkin type
$R=A_{16}+A_{2}+A_{1}$.
We put
$$
f\;:=\;17\, t^3 -18\, t^2 -228\, t+536,
$$
which has two non-real roots $\alpha, \bar{\alpha}$ and a real root $\beta$.
In~\cite{MR1900779},
it is shown that
$\MMM(R)$ consists of three  connected components
$\MMM(R)_\alpha$, $\MMM(R)_{\bar{\alpha}}$ and $\MMM(R)_\beta$,
and that there exists a homogeneous polynomial $\Phi(x_0, x_1, x_2)$
of degree $6$ with coefficients in $F_f$ such that
the conjugate sextics 
$$
C_\alpha=\{\Phi^\alpha=0\},
\quad
C_{\bar{\alpha}}=\{\Phi^{\bar\alpha}=0\},
\quad
C_\beta=\{\Phi^\beta=0\}
$$
are  members of $\MMM(R)_\alpha$, $\MMM(R)_{\bar{\alpha}}$ and $\MMM(R)_\beta$,
respectively.
On the other hand,
by the method described in the previous section, 
we calculate  that
$\overline{\Ms}\spsharp=\{[M^0]\}$
and $\Ns(M^0)=\{N^1, N^2\}$, where 
$$
N^1=\Lat  [10,4,22]\quand
N^2=\Lat  [6,0,34].
$$
The set  $F([M^0], N^1)$ consists of two non-real orbits,
while the set $F([M^0], N^2)$ consists of a single real orbit.
Since the complex conjugation induces a homeomorphism
$(\Pt, C_\alpha)\cong (\Pt, C_{\bar{\alpha}})$,
the invariants $N_{C_\alpha}$ and $N_{C_{\bar\alpha}}$ 
must be equal.
Hence  
$N_{C_\alpha}$ and $N_{C_{\bar\alpha}}$ 
are isomorphic to $N^1$,
while  $N_{C_\beta}$ is isomorphic to $N^2$.
Since $N^1$ and $N^2$ are not isomorphic, 
we conclude that   $(C_\alpha, C_\beta)$
is an arithmetic Zariski pair.
\end{example}
\begin{example}
Consider the Dynkin type
$R=A_{16}+A_{3}$.
In~\cite{MR1900779},
it is shown that
$\MMM(R)$ consists of two  connected components
$\MMM\sb +$ and $\MMM\sb -$,
and that there exist members $C\sb +$ of $\MMM\sb +$ and $C\sb -$ of $\MMM\sb -$
that are conjugate in $\Q(\sqrt{17})$.
On the other hand,
we calculate  that
$\overline{\Ms}\spsharp=\{[M^0]\}$
and   $\Ns(M^0)=\{N^1, N^2\}$, where 
$$
N^1=\Lat  [4,0,34]\quand
N^2=\Lat  [2,0,68].
$$
Each of  $F([M^0], N^1)$ and $F([M^0], N^2)$
consists of a single real orbit.
Therefore $(C\sb +, C\sb-)$
is an arithmetic Zariski pair.
\end{example}
\begin{example}
Suppose that  $R=A_{18}+A_{1}$.
We put 
$$
f\;:=\;19\, t^3+50\, t^2 +36\, t +8,
$$
which has two non-real roots $\alpha, \bar{\alpha}$ and a real root $\beta$.
Again by~\cite{MR1900779},
the moduli space $\MMM(R)$
consists of three  connected components
$\MMM(R)_\alpha$, $\MMM(R)_{\bar{\alpha}}$, $\MMM(R)_\beta$
that have members
$C_\alpha=\{\Phi^\alpha=0\}$,
$C_{\bar{\alpha}}=\{\Phi^{\bar{\alpha}}=0\}$,
$C_\beta=\{\Phi^\beta=0\}$, respectively,  
for some homogeneous polynomial $\Phi$ with coefficients  in $F_f$.
On the other hand,
we have 
$\overline{\Ms}\spsharp=\{[M^0]\}$
and $\Ns(M^0)=\{N^1, N^2\}$,
where
$$
N^1=\Lat  [8,2,10]\quand
N^2=\Lat  [2,0,38].
$$
The set  $F([M^0], N^1)$ consists of two non-real orbits,
while the set $F([M^0], N^2)$ consists of a single real orbit.
Hence  
$(C_\alpha, C_\beta)$
is an arithmetic Zariski pair.
\end{example}
For the cases $R=A_{15}+A_{4}$ and $R=A_{19}$
that are also treated in~\cite[Theorem~5.8]{MR1900779},
the situation is as follows.
\begin{example}
Suppose that $R=A_{15}+A_{4}$.
We have
$\overline{\Ms}\spsharp=\{[M^0], [M^1]\}$,
where $M^1$ is an overlattice of $M^0$ with index $2$.
We have  
$\Ns([M^0])=\{N^0\}$
with  $N^0=\Lat  [8,4,22]$, and   
$F([M^0], N^0)$ consists of two non-real orbits,
while we have 
$\Ns([M^1])=\{N^1\}$
with  $N^1=\Lat  [2,0,20]$, and  
$F([M^1], N^1)$ consists of a single  real orbit.
According to Yang's list~\cite{MR1387816},
there exist two  configurations
of maximizing sextics of type $A_{15}+A_{4}$.
By~\cite{MR1900779},
there exist 
 members $C$ and $\overline{C}$ 
 of distinct connected components of $\MMM(R)$ 
 that are conjugate in $\Q(\sqrt{-1})$.
The complex conjugation yields a homeomorphism
  $(\Pt, C)\cong (\Pt, \overline{C})$.
 Hence 
 we must have $N_{C}\cong N_{\overline{C}}\cong N^0$.
\end{example}
\begin{example}
Suppose that $R=A_{19}$.
We have
$\overline{\Ms}\spsharp=\{[M^0]\}$
and $\Ns([M^0])=\{N^0\}$
with   $N^0=\Lat  [2,0,20]$.
The set $F([M^0], N^0)$ consists of two real orbits.
According to~\cite{MR1900779},
there exist
 members $C\sb +$ and $ C\sb -$ of  $\MMM(R)$
 belonging to
 the distinct connected components
  that 
are conjugate in $\Q(\sqrt{5})$.
Our invariant fails to distinguish $(\Pt, C_+)$ and $(\Pt, C_-)$ topologically, because we have 
$N_{C\sb+}\cong N_{C\sb-}\cong N^0$.
It would be an interesting problem to determine whether
$(\Pt, C\sb+)$ and $(\Pt, C\sb-)$ are homeomorphic or not.
\end{example}
\section{A singular $K3$ surface defined over  a number field}\label{sec:singK3}
Let $Y$ be a complex $K3$ surface or a complex abelian surface
such that the transcendental lattice $T(Y)$  is
of rank $2$.
Then $T(Y)$ is an even positive-definite lattice.
Moreover the Hodge structure
$$
T(Y)\otimes \C=H^{2, 0}(Y)\oplus\overline{H^{2, 0}(Y)}
$$
of $T(Y)$ defines a canonical orientation on $T(Y)$;
namely, an ordered basis $e_1, e_2$ of $T(Y)$ is positive
if the imaginary part of the complex number $(e_1, \omega_Y)/(e_2, \omega_Y)$
is positive,
where $\omega_Y\in H^{2, 0}(Y)$
is the cohomology class of a non-zero holomorphic $2$-form of $Y$.
We denote by $\wt{T}(Y)$ the \emph{oriented transcendental lattice} of $Y$.
\begin{definition}
A (smooth) $K3$ surface $X$ defined over a field  $k$ of
 characteristic $0$ is called \emph{singular}
if the Picard number of $X\otimes \bar{k}$ is $20$.
\end{definition}
If $X$ is a \emph{complex} singular $K3$ surface,
then we have the oriented transcendental lattice $\wt{T}(X)$.
We have the following important theorem due to Shioda and Inose~\cite{MR0441982}:
\begin{theorem}\label{thm:SI}
The correspondence $X\mapsto \wt{T}(X)$
yields a bijection from the set of isomorphism classes of complex singular $K3$ surfaces
to the set of isomorphism classes of  even positive-definite oriented lattices of rank $2$.
\end{theorem}
Notice that, if $C$ is a complex maximizing sextic,
then the minimal resolution $X_C$ of the double covering $Y_C\to \Pt$
that branches  exactly along $C$ is a complex singular $K3$ surface,
and $T(X_C)$ is isomorphic to $N_C$.
\par
\medskip
Let $X$ be a singular $K3$ surface defined over a number field $F$.
For an embedding $\sigma$ of $F$ into  $\C$,
we denote by $X\sp\sigma$ the complex $K3$ surface obtained from $X$ by
$\sigma$.
\par
\medskip
The following is a special case of~\cite[Theorem~3]{tssK3}.
\begin{proposition}\label{prop:55}
There exist a singular $K3$ surface $X$ defined over
a number field $F$ and two embeddings $\tau$ and $\tau\sprime$ of $F$ into $\C$ such that
$$
\wt{T}(X\sp\tau)\;\cong\;\oriLat  [2,1,28]
\quand
\wt{T}(X\sp{\tau\sprime})\;\cong\; \oriLat  [8,3,8].
$$
\end{proposition}
\begin{proof}
We put
$$
K:=\Q(\sqrt{-55})\;\;\subset\;\;\C,
$$
and denote by $\Z_K$ the ring of integers of $K$.
For a number field $L$ containing $K$,
we denote by $\Emb (L/K)$ the set of embeddings of $L$ into $\C$
whose  restrictions to $K$ are the identity of $K$.
We define fractional ideals $I_0$, \dots, $I_3$ of $\Z_K$ by the following:
$$
\renewcommand{\arraystretch}{1.2}
\begin{array}{cclcl}
I_0 &:=& \Z_K = \Z+\Z\tau_0&\textrm{where} &\tau_0:=(1+\sqrt{-55})/2, \\
I_1 &:= &\Z+\Z\tau_1 &\textrm{where} &\tau_1:=(3+\sqrt{-55})/4, \\
I_2 &:=&\Z+\Z\tau_2&\textrm{where} &\tau_2:=(5+\sqrt{-55})/8, \\
I_3&:= &\Z+\Z\tau_3&\textrm{where} &\tau_3:=(1+\sqrt{-55})/4.
\end{array}
$$
The ideal class group $\Cl_K$ of  $\Z_K$ 
is a cyclic group of order $4$
generated by the  class $[I_1]$,
and we have
$[I_2]=[I_1]^2$ and  $[I_3]=[I_1]^3$.
We consider the Hilbert class polynomial
\begin{eqnarray*}
\HHH (t)&:=&(t-j(\tau_0))(t-j(\tau_1))(t-j(\tau_2))(t-j(\tau_3))\\
&=& {t}^{4}+13136684625\,{t}^{3}-20948398473375\,{t}^{2}+\mystruth{14pt}\\&&\qquad
 +172576736359017890625\,t-18577989025032784359375
\end{eqnarray*}
of $\Z_K$, and the Hilbert class field 
$$
H:=K[t]/(\HHH(t))
$$
of $K$ (see Cox~\cite{MR1028322}).
We put
$$
\gamma:= t\bmod (\HHH)\in H.
$$
Then we have
$\Emb (H/K)=\{\sigma_0, \sigma_1, \sigma_2, \sigma_3\}$,
where $\sigma_i$ is the embedding  defined by
$\sigma_i(\gamma)=j(\tau_i)\in \C$.
Consider the elliptic curve 
\begin{equation}\label{eq:E}
E\;\;:\;\; y^2+xy=x^3-\frac{36}{\gamma-1728} x-\frac{1}{\gamma-1728}
\end{equation}
defined over $H$
(see Silverman~\cite[page 52]{MR817210}).
Then we have $j(E)=\gamma\in H$, and
hence  $j(E^{\sigma_i})=j(\tau_i)$ holds
for $i=0, \dots, 3$,
where $E^{\sigma_i}$ is the complex elliptic curve
defined by~\eqref{eq:E} with $\gamma$ replaced by $\sigma_i(\gamma)=j(\tau_i)$.
Therefore we have an isomorphism of Riemann surfaces
\begin{equation}\label{eq:Esigma}
E^{\sigma_i}\;\;\cong\;\;\C/ I_i
\end{equation}
for $i=0, \dots, 3$. 
We then put
$$
A:=E\times E.
$$
Note that $T(A)$ is of rank $2$.
By means of a double covering of the Kummer surface $\Km(A)$ of $A$, 
Shioda and Inose~\cite{MR0441982} constructed a singular $K3$ surface $X$
defined over a finite extension $F$ of $H$ with the following properties
(see also~\cite[Propositions 6.1 and 6.4]{tssK3}).
$$
\parbox{10cm}{
For any  $\sigma\in \Emb(F/K)$,
the oriented trancendental lattice $\wt{T}(X\sp\sigma)$
is isomorphic to the oriented transcendental lattice $\wt{T}(A\sp\sigma)$
of the complex abelian surface $A\sp\sigma=E\sp\sigma\times E\sp\sigma$.}
$$
See Inose~\cite{MR578868} and Shioda~\cite{Shioda_murre}
for an explicit defining equation of $X$.
\par
\medskip
The oriented lattice $\wt{T}(A\sp\sigma)$
is calculated by Shioda and Mitani~\cite{MR0382289}.
Suppose that the restriction of $\sigma\in \Emb(F/K)$ to $H$ is $\sigma_i$.
Then we have
$$
A\sp\sigma\;\;\cong\;\; \C/I_i\times \C/I_i\;\; \cong \;\; \C/I_i^2\times \C/I_0\;\;
\cong\begin{cases}
\C/I_0\times \C/I_0 &\textrm{if $i=0$ or $i=2$,}\\
\C/I_2\times \C/I_0 &\textrm{if $i=1$ or $i=3$,}
\end{cases}
$$
by~\eqref{eq:Esigma} and~\cite[(4.14)]{MR0382289}.
Hence, by~\cite[Section 3]{MR0382289}, we have
$$
\wt{T}(A\sp\sigma)\cong\begin{cases}
\oriLat [2,1,28] &\textrm{if $\sigma|H$ is $\sigma_0\;\textrm{or}\;\sigma_2$,}\\
\oriLat [8,3,8] &\textrm{if $\sigma|H$ is $\sigma_1\;\textrm{or}\;\sigma_3$,}
\end{cases}
$$
(see also~\cite[\S6.3]{tssK3}.)
Thus we obtain the hoped-for $X$ and $\tau$, $\tau\sprime$.
\end{proof}
\begin{remark}\label{rem:orirev}
Note that the orientation reversing does not change  the isomorphism classes
of the oriented lattices $\oriLat[2,1,28]$ and $\oriLat[8,3,8]$
(see Remark~\ref{rem:forget}).
Hence, by
Theorem~\ref{thm:SI}, 
if a complex singular $K3$ surface $Y$
satisfies  $T(Y)\cong \Lat[2,1,28]$
(resp.~$T(Y)\cong \Lat[8,3,8]$),
then $Y$ is isomorphic to the complex $K3$ surface $X\sp\tau$
(resp.~to the complex $K3$ surface $X\sp{\tau\sprime}$) in Proposition~\ref{prop:55}
\end{remark}
Using Proposition~\ref{prop:55} and Remark~\ref{rem:orirev},
we obtain the following example of arithmetic Zariski pairs.
\begin{example}
Consider the Dynkin type $R=A_{10}+A_{9}$.
We have
$$
\overline{\Ms}\spsharp=\{[M^0], [M^1]\},
$$ 
where $M^1$ is an overlattice of $M^0$ with index $2$.
We then have
\begin{eqnarray*}
\Ns ([M^0])&=&\{\;\Lat  [10, 0, 22],\;\; \Lat  [2,0,110]\;\}\quand \\
\Ns ([M^1])&=&\{\;\Lat  [2,1,28],\;\; \Lat  [8,3,8]\;\}, 
\end{eqnarray*}
and each of the sets 
$$
\renewcommand{\arraystretch}{1.2}
\begin{array}{lll} 
&F([M^0], \Lat  [10, 0, 22]),
\quad
&F([M^0], \Lat  [2,0,110]),\\
&F([M^1], \Lat  [2,1,28]),
\quad
&F([M^1], \Lat  [8,3,8])
\end{array}
$$
consists of a single real orbit.
In particular, the number of the  connected components
of $\MMM(R)$ is four.
Let $C$ and $C\sprime$ be members of 
the connected components of $\MMM(R)$ corresponding to 
$F([M^1], \Lat  [2,1,28])$ and 
$F([M^1], \Lat  [8,3,8])$,
respectively.
Note that we have 
$$
N_C\;\cong\; T(X_C)\;\cong\; \Lat  [2,1,28]
\quand
N_{C\sprime}\;\cong\; T(X_{C\sprime})\;\cong\; \Lat  [8,3,8].
$$
By Remark~\ref{rem:orirev},
we see that $X_C$ is isomorphic to  $X\sp{\tau}$
and $X_{C\sprime}$ is isomorphic to $X\sp{\tau\sprime}$.
Consider the composites
$$
\phi_C\;:\; X_C\maprightsp{}  Y_C\maprightsp{\pi_C} \Pt
\quand
\phi_{C\sprime}\;:\; X_{C\sprime}\maprightsp{}  Y_{C\sprime}\maprightsp{\pi_{C\sprime}} \Pt
$$
of the finite double  coverings branching along $C$ and $C\sprime$
and the minimal desingularizations.
Since $X_C\cong X\sp{\tau}$, there exists a morphism
$\phi_L: X\otimes L\to \Pt$
with the Stein factorization
$$
\phi_L\;:\; X\otimes L \maprightsp{}  Y_L \;\maprightsp{\pi_L}\;\Pt
$$
defined over a finite extension $L$ of $F$ such that, 
for some embedding $\theta$ of $L$ into $\C$
satisfying $\theta|F=\tau$,
the morphism
$$
\phi_L^\theta\;:\;
(X\otimes L)^\theta=X\sp\tau  
\maprightsp{}  Y_L^\theta 
\;\maprightsp{\pi_L^\theta} \;
\Pt
$$
is isomorphic to $\phi_C$.
In particular,  the branch curve $B$
of the finite double  covering $\pi_L^\theta$
is isomorphic to $C$ as  a complex plane curve.
Let $\theta\sprime$ be an embedding of $L$ into $\C$
such that $\theta\sprime|F=\tau\sprime$,
and consider the morphism 
$$
\phi_L^{\theta\sprime}\;:\;
(X\otimes L)^{\theta\sprime}=X\sp{\tau\sprime}  \maprightsp{}  
Y_L^{\theta\sprime}\;\maprightsp{\pi_L^{\theta\sprime}} \;
\Pt.
$$
Since the branch curve $B\sprime$ of $\pi_L^{\theta\sprime}$ is conjugate to
the branch curve $B$ of $\pi_L^\theta$,
it is a maximizing sextic of type $A_{10}+A_{9}$.
Since $ (X\otimes L)^{\theta\sprime}=X\sp{\tau\sprime}$  is isomorphic to $X_{C\sprime}$, 
the morphism
$\phi^{\theta\sprime}$
must be isomorphic to $\phi_{C\sprime}$, and hence 
$B\sprime$ is isomorphic to  $C\sprime$
as a complex plane curve.
Therefore 
the conjugate pair $(B, B\sprime)$
of plane curves 
 is isomorphic to
the pair 
$(C, C\sprime)$ with $N_C\not\cong N_{C\sprime}$, 
and thus yields an example of arithmetic Zariski pairs.
\end{example}
\bibliographystyle{plain}
\def\cprime{$'$} \def\cprime{$'$} \def\cprime{$'$}

\end{document}